\documentstyle[11pt]{article}

\newtheorem{th}{Theorem}[section]
\newtheorem{de}[th]{Definition}

\newtheorem{lem}[th]{Lemma}
\newtheorem{co}[th]{Corollary}
\newtheorem{re}[th]{Remark}

{

\begin{document}

\large
\title {\bf THE PENALTY METHOD FOR VARIATIONAL INEQUALITIES WITH
NONSMOOTH UNBOUNDED OPERATORS IN BANACH SPACE}
\author{\bf Ya. I. Alber $^*$\\
Department of Mathematics\\
Technion-Israel Institute of Technology\\
Haifa 32000, Israel}
\date{}
\maketitle

{\bf Abstract} -
The existence of a solution, convergence and stability of the  penalty
method for variational inequalities with nonsmooth unbounded
uniformly  and
properly monotone operators in Banach spase $B$
are  investigated.   All  the  objects  of  the
inequality  -  the  operator  A,  "the  right-hand  part" $f$  and  the  set
of constrains $\Omega $ - are  to  be  perturbed. The stability  theorems
are formulated in terms of geometric characteristics of the  spaces
$B$ and $B^*$.
The results of this paper are continuity and generalization of
the Lions' ones, published  earlier in  \cite{l}. They are new even in
Hilbert spaces.

\section{Preliminaries}
\setcounter{equation}{0}

The goal of our paper is to prove the  existence of a solution,
convergence and stability
of the penalty method for the variational inequalities (\ref{p1}) with the
{\bf unbounded }
operator $A$ in Banach space $B$. At the same time, formulation of our
problem includes other complicated  aspects: the operator $A$ is able to be
properly monotone (i.e. degenerate) and
nonsmooth (possibly, multivalued and even discontinuous), and its domain
does not necessarily coincide with the whole space $B.$
A quite new moment here is the investigation of this
method  on approximately given sets $\Omega.$ In addition,
the operator $A$ and "the right-hand part" $f$ are to  be perturbed.
Apparently, it is
the most general formulation of a problem about
the penalty method. Earlier,
only regularization method has been considered under these
conditions (see \cite{a1,an3}). All the results are new even in Hilbert
spaces.

Let $B$ be  a  real  reflexive  Banach  space  with  a  strictly  convex
conjugate space\\
------------------------\\
$^{*}$ This research was supported in part by the Ministry of Science Grant
3481-1-91 and by the Ministry of Absorption Center for Absorption in Science.\\
$B^*$ and  $<y,x>$ a dual product in $B$, i.e.a pairing
between $y \in B^*$ and $ x \in B.$
Let the signs $||\cdot||$ and $ \;||\cdot||_{B^{*}}$  denote
the norms in the Banach spaces $B$ and $B^{*}$, respectively.
Suppose that $A $ is a monotone
operator from $B$ to $B^*$ with domain  $D(A) \subseteq B,\; \Omega $  is  a
closed convex set, $\Omega  \subseteq  $ {int}$D(A)$,  {int}$\Omega  \neq
\emptyset.$

\begin {de}  \label{gpo}
An element $x^* \in  \Omega $  is  said  to  be  a  solution  of  the
variational inequality
\begin{equation} \label{p1}
<Ax-f,y-x> \ge  0, \;\;\;\forall y \in  \Omega , \;\;f \in  B^*, \;\;x
\in D(A)
\end{equation}
if there exists $z \in \bar A x^*$ such that
\begin{equation} \label{p2}
<z-f,x-x^*> \ge  0,\;\; \forall x \in  \Omega
\end{equation}
where $\bar A$  is the maximal monotone expansion of $A$ in $D(A)$. \\
\end {de}

This definition does not depend on the smoothness of  the  operator  $A $
because for a monotone operator there always exists  its  maximal  monotone
extension (according to Zorn's Lemma).

Note that (\ref {p2}) is equivalent to the inequality
\begin{equation} \label{p25}
<z-f,x-x^*> \ge  0,\;\; \forall x \in  \Omega, \;\; \forall z \in Ax.
\end{equation}

Recall that the operator $A$ is said to be\\
(i) $\;$monotone if $<y_{1}-y_{2},x_{1}-x_{2}> \ge  0, \forall y_{1} \in
Ax_{1},
\forall y_{2} \in  Ax_{2}, \forall x_{1}, x_{2} \in  D(A);$ \\
(ii) $\;$strictly  monotone  if it is monotone and
$<y_{1}-y_{2},x_{1}-x_{2}> > 0, \forall y_{1}
\in Ax_{1}, \forall y_{2} \in Ax_{2},
\forall x_{1}, x_{2} \in  D(A), x_{1} \neq  x_{2}$;\\
(iii) proper  monotone  (strictly  proper  monotone)  if  it  is   monotone
(strictly monotone) and there is no intensification of conditions (i)
and (ii), respectively;\\
(iv) $\;$uniformly  monotone  if $<y_{1}-y_{2},x_{1}-x_{2}> \ge  \psi
(||x_{1}-x_{2}||),\; \forall y_{1} \in  Ax_{1},
\forall y_{2} \in  Ax_{2}, \forall x_{1}, x_{2} \in  D(A)$, where $\psi (t)$
is  a  continuous  positive function and $\psi (0) = 0$; \\
(v) $\;$coercive in $D(A)$ if $\lim <Ax,x> / ||x||  = \infty $  as $||x||
\rightarrow \infty ,\; x \in  D(A)$;\\
(vi) $\;$$\Omega $-coercive in $D(A)$  if  there  exists  a  point $x_{0} \in
\Omega $  such  that
$$\lim \frac {<Ax,x-x_{0}>}{||x||}  = \infty $$
as $||x|| \rightarrow \infty ,
x \in  D(A)$; \\
(vii)  maximal monotone  if its graph $G(A)$ is a maximal monotone set;\\
(viii) bounded if it carries bounded sets of $D(A)$ to bounded sets of $B^*$.\\
\\

The penalty method replaces the problem  of  finding  the  variational
inequality solution on the set $\Omega  \subset  D(A)$ with the problem of
finding solutions to the equations on  the  whole  domain $D(A)$,  by
introducing  a  penalty operator with an arbitrarily large parameter.
The penalty term is equal to zero on the set $\Omega $ and tends to
infinity at each point $x \notin \Omega $ as the  penalty
parameter increases indefinitely. \\
\\
We use a penalty operator of the kind
$$\mu (x) = J(x-P_{\Omega }x)$$
which has been proposed earlier by Lions in \cite{l}. Here $x \in  B,
\;P_{\Omeg
\rightarrow \Omega $ is a  metric projection
operator on the set $\Omega $ which yields the correspondence between each
point $x \in  B$  and
its nearest one $\bar x$  in $\Omega$  (best approximation) according
to minimization problem
\begin{equation} \label{p3}
||x - \bar x||  = ||x - P_{\Omega }x|| = \inf ||x-\xi|| , \;\;\forall \xi
\in  \Omega,
\end{equation}
$J: B \rightarrow  B^*$ is normalized duality mapping, determined by the
relations
$$<Jx,x> = ||Jx||_{B^*}||x||  = ||x||^{2}.$$
Under our conditions the projection operator exists in every point $x \in  B$
and it is single-valued.  Note  also  that  a  duality  mapping  exists  in
arbitrary Banach space and is most simply described in spaces $l ^{p},\;L^{p},
\;W^{p}_{m},\; \infty > p > 1$:
\begin{itemize}
\item  (i) $ l^p:  Jx = ||x||^{2-p}_{l^p} y\in l^q, {\quad } x =
\lbrace  x_1, x_2,...\rbrace  $,
 $ y =   \lbrace  x_1 ||x_1||^{p-2}, x_2 ||x_2||^{p-2},...\rbrace, \\
 p^{-1} + q^{-1} = 1 ,$
\item  (ii) $ L^p:  Jx = ||x||^{2-p}_{L^p} |x|^{p-2}x \in L^q $
 \item  (iii) $ W^p_m:  Jx = ||x||^{2-p}_{W^p_m} \sum (-1)^ {|\alpha
|}D^\alpha
   (|D^\alpha x |^{p-2} D^\alpha x ) \in W^q_{-m} .$
 \end{itemize}
In Hilbert space  $J$ is  an identity operator.

The dual mapping has a lot of remarkable properties.  It is a monotone
and  bounded  operator  in  arbitrary  Banach  space,  a   continuous   and
single-valued operator in  smooth  Banach  space,  a  uniformly  continuous
operator on each bounded set  in  uniformly  smooth  Banach  space,  and  a
uniformly monotone operator on each bounded set in uniformly  convex  Banach
space.  These properties were described long time ago in  plenty of books
concerning
the theory of monotone  operators,  but  the  quantitative  description  of
uniform properties mentioned above has been obtained only  recently
\cite{an1,an2}.
In particular, if $\delta _{B}(\epsilon )$ and $\rho _{B}(\tau)$ are modulus
of convexity and  modulus  of
smoothness of Banach space $B$ (respectively) and if $||x|| \leq  R$
and $||y|| \leq  R$, then
$$ <Jx-Jy,x-y>) \ge  (2L)^{-1}\delta _{B}(||x-y||C_{2}),\; 1 < L < 3.18,\;
C_{2} = {max} \lbrace 1,R \rbrace, $$
\begin{equation} \label{p18}
||Jx-Jy||_{B^*} \leq  C_{2}g^{-1}_{B^*}(2LC_{2} ||x-y||),\; g_{B}(\epsilon ) =
\delta _{B}(\epsilon )/\epsilon,
\end{equation}
\begin{equation} \label{p19}
<Jx-Jy,x-y> \leq  8 ||x-y||^{2} + 8C_{1} \rho _{B}(||x-y||),\; C_{1} =
{max}\lbrace L,R \rbrace ,
\end{equation}
where $g^{-1}_{B}(\cdot )$ is an inverse function of $g_{B}(\epsilon )$.  The
dual estimates are valid too \cite{an1,an2}.

Let us mark the following properties \cite{d,an4,p}:\\
(i) the function $\delta _{B}(\epsilon )$ is defined on the interval  [0, 2],
continuous and  increasing in a uniformly convex Banach space,
$\delta _{B}(0) = 0,\;\delta _{B}(\epsilon) <1 ,$ \\
(ii) the function $\rho _{B}(\tau)$ is defined on the interval [0, $\infty$),
convex, continuous and nondecreasing, $\rho _{B}(0) = 0,$\\
(iii) the function $g_{B}(\epsilon )$ is continuous and increasing on (0,2],
$g_{B}(0) = 0$.

In addition, one can  consider that $c_1 \epsilon^{\gamma}
\leq  g_{B}(\epsilon )
\leq c_2 \epsilon ,\; \gamma \ge 1$ (the upper estimate follows from the
inequality $ \delta _{B}(\epsilon) \leq \delta _{H}(\epsilon) \leq
\epsilon ^2 /4,$ where $H$ is Hilbert space; the lower estimate follows
from  \cite{fi,p}).

Let us show that estimate
\begin{equation} \label{p4}
<Jx-Jy,x-y> \leq {\bar C}_{1} \rho _{B}(||x-y||)
\end{equation}
occurs also. Indeed  \cite{d},  in  Hilbert  space $H,\; \rho _{H}(\tau) =
(1 + \tau^{2})^{1/2}-1.$
Therefore,
$$\rho _{B}(\tau) \ge  \rho _{H}(\tau) \ge ((1+4R^{2})^{1/2} +
1)^{-1} \tau^{2}, $$
because $\tau \leq  2R.$  Thus (\ref{p4}) is
satisfied from (\ref{p19}) with constant
$${\bar C}_{1} = 8 ((1+4R^{2})^{1/2}+1) + {max} \lbrace L,R \rbrace .$$

The operator $\mu (x)$ is monotone \cite{l}, i.e.,
$$<J(x-P_{\Omega }x) - J(y-P_{\Omega }y), x-y> \ge  0,$$
and, besides,
$$<J(x-P_{\Omega }x) - J(y-P_{\Omega }y), \bar x - \bar y> \ge  0,$$
\begin{equation} \label{p5}
<J(x-P_{\Omega }x), P_{\Omega }x-\xi> \ge  0,\;\; \forall \xi  \in  \Omega .
\end{equation}
We use the following statements on the theorems below.

\begin{lem}\label{l1}
$\bar x \in  \Omega $ is a projection of the point $x \in  B$ on $\Omega $
if and only if the inequality
\begin{equation} \label{p6}
||x - \bar x||^{2} \leq  <J(x - \bar x),x-\xi>,\;\; \forall \xi  \in  \Omega
\end{equation}
is satisfied.
\end{lem}

{\bf Proof}.  Indeed, from (\ref{p6}) it follows immediately that
$$||x - \bar x|| \leq ||x - \bar x||^{-1} <J(x - \bar x),x-\xi>
 \leq ||x - \xi||,\;\; \forall \xi  \in  \Omega, $$
i.e., (\ref{p3}) is valid.  Therefore $ \bar x  = P_{\Omega }x$.  Inversely,
if $ \bar x  = P_{\Omega }x$,  then  by virtue of (\ref{p5}) we have
$$<J(x - \bar x), \bar x - \xi> = <J(x - \bar x), \bar x - x> + $$
$$<J(x - \bar x), x - \xi> = - ||x - \bar x||^{2} + <J(x - \bar x), x - \xi>
\ge  0$$
and (\ref{p6}) occurs.

\begin{re}
Let $x \in B$ be a fixed and $\xi  \in  \Omega $ an arbitrary point,
$$P(\xi ) = <J(x - \bar x), x - \xi>,\;\;\; Q(\xi ) = <J(x - \bar x),
\bar x- \xi> .$$
Then $P(\xi ) \ge  0,\; Q(\xi ) \ge  0$  and
$P(\xi ) - Q(\xi ) = d = const, d = ||x - \bar x||^{2}$.
\end{re}

\begin{lem}\label{l2}
\cite{an3,a1} Let $B$  be an arbitrary Banach space, $F: D(F) \subseteq  B
\rightarrow  B^* $ be a monotone operator, $x^{0} \in {int} D(F)$.
Then  there  exist  the  constant $r_{0} > 0$ and the closed ball
$S(r_{0},x^{0}) \subset  D(F)$ such that for all $x \in  D(F)$  the estimate
$$<F(x),x-x^{0}> \ge  r_{0}||F(x)|| - c_{0} (||x-x^{0}|| + r_{0}), \;\;
c_{0} = \sup ||F(\xi)|| < \infty \;, \xi  \in  S(r_{0},x^{0}) $$
is satisfied.
\end{lem}

\begin{lem}\label{l3}
\cite{a1,an4} Let $B$ be a uniformly convex Banach space, $\Omega _{1}$ and
$\Omega _{2}$  be
convex closed sets, $x \in  B,\; {\cal H}(\Omega _{1},\Omega _{2}) \leq
\sigma $, where ${\cal H}(\Omega _{1},\Omega _{2})$ is  the  Hausdorff
distance between $\Omega _{1}$ and $\Omega _{2}$.  If $\delta _{B}(\epsilon)$
is a modulus of convexity of space $B$,
and $\delta ^{-1}_{B}(\cdot )$ is it's inverse function, then
$$||P_{\Omega _{1}}x - P_{\Omega _{2}}x|| \leq  C\delta ^{-1}_{B}
(4LC_{1}\sigma), 1 < L < 3.18 ,$$
$$C = 2 \max  \lbrace 1,||x - \bar x_{1}||,||x - \bar x_{2}|| \rbrace,\;\;
\bar x_{1} = P_{\Omega _{1}}x,\; \bar x_{2} = P_{\Omega _{2}}x, $$
$$C_{1} = 2 \max \lbrace ||x - \bar x_{1}||,||x - \bar x_{2}|| \rbrace. $$
\end{lem}

\section{Variational  inequalities on the exactly given sets.  Convergence
of the penalty method}
\setcounter{equation}{0}

First consider the penalty method for variational  inequality with maximal
monotone operator $A$ on  the exactly given set $\Omega $:
\begin{equation} \label{p7}
Ax + \epsilon ^{-1}J(x-P_{\Omega }x) = f,\;\; x \in  D(A),\; \epsilon > 0,\;
\epsilon \rightarrow 0 .
\end{equation}

\begin {de}  \label{gpo}
The element $x_{\epsilon } \in  D(A)$ is said to be  a  generalized  solution
to  equation (\ref{p7})
for fixed $\epsilon  > 0$ if there exists $z_{\epsilon } \in  A x_
{\epsilon }$ such that
\begin{equation} \label{p8}
<z_{\epsilon } - f, x - x_{\epsilon }> + \epsilon ^{-1}<J(x_{\epsilon } -
P_{\Omega }x_{\epsilon }), x - x_{\epsilon }> \ge  0,\;\; \forall x \in  D(A).
\end{equation}
\end {de}

Next we will obtain the existence theorem for the variational inequality
(\ref{p1}) with  nonsmooth unbounded operator $A$ by means of the penalty
sequence $x_{\epsilon }.$  Note that the most important and most difficult
moment in the proof of convergence and stability of any approximation
method is to establish uniform boundedness of the operator $A$ of
variational inequality on corresponding approximation solutions. Emphasize
that the technique developed before in \cite{l,glt} for (\ref{p1})  and
(\ref{p7}) turned out to be
unsuitable in considered (unbounded) case especially if the set of
constrains is perturbed. To solve this problem we propose to use Lemma
\ref{l2}. By means of this Lemma one can show that a priori unbounded
operator $A$ is {\it always uniformly bounded \/}
on the solutions of penalty equation (\ref{p7}). The proving gives
constructive method of approximation of the solution $x^*.$

\begin {th}\label{t1}
$\;$Let  $A$ be a strictly monotone operator,$\;$ $\Omega$  a bounded  set,
$\;$
$\theta \in int D(A)\; (\theta $ is an origin of the space $B$). Then  there
exists  a  unique
solution $x^*$ of the variational inequality (\ref{p1}), and generalized
solutions $x_{\epsilon }$
of the equation (\ref{p7}), which are uniformly bounded, converge weakly
to $x^*$
as $\epsilon \rightarrow 0$, and  $||x_{\epsilon }-P_{\Omega }x_{\epsilon }||
= O(\epsilon )$.
\end {th}

{\bf Proof.} The operator
$$Tx =  Ax +  \epsilon ^{-1}J(x-P_{\Omega }x) $$
is coercive and strictly monotone in $D(A)$. Indeed, let $x \in D(A).$ Then
$$<Tx,x> = <Ax,x> + \epsilon ^{-1} < J(x-P_{\Omega }x),x> = <Ax -
A\theta,x>$$
$$ +  <A\theta,x> + \epsilon ^{-1} < J(x-P_{\Omega }x),x - P_{\Omega }x> +
\epsilon ^{-1} < J(x-P_{\Omega }x), P_{\Omega }x> $$
$$\ge - ||A\theta||_{B^*}||x|| + \epsilon ^{-1}||x - P_{\Omega }x ||^2 -
\epsilon ^{-1}||x - P_{\Omega }x ||||P_{\Omega }x||.$$
The coerciveness of $T$ follows from a local boundedness of maximal monotone
operator $A$ at each point of $\;intD(A)$, from a boundedness of
the set $\Omega$ and from the condition $\theta \in int D(A).$
A strict monotonicity of $T$ follows from a monotonicity of the penalty
operator $\mu (x)$ and strictly monotonicity of $A.$
It is obviously that
$$D(A) \bigcap D(G) \neq \emptyset,$$
where $G(x) = J(x-P_{\Omega }x), $
and by virtue of \cite{rc}, the operator $T$ is maximal monotone.

The existence and uniqueness of the generalized
solution $x_{\epsilon }$ of the penalty equation (\ref{p7}) for
fixed $\epsilon > 0$  follow  from  the strict  monotonicity
and  coerciveness of the maximal monotone operator $T$, from the
reflexiveness of $B$ and strict  convexity of $B^{*}$,
and from the condition $\theta \in D(A)$ (see \cite{b,ar}).

By virtue of Lemma \ref{l1},
if $\xi  \in  \Omega ,\; z_{\epsilon } \in A x_{\epsilon }, $
we can write for the generalized solution
$x_{\epsilon }$
\begin{equation} \label{p9}
||x_{\epsilon } - \bar x_{\epsilon }||^{2} \leq  <J(x_{\epsilon } -
\bar x_{\epsilon }), x_{\epsilon } - \xi>\; \leq \; \epsilon <f-z_{\epsilon },
x_{\epsilon } - \xi>.
\end{equation}
Let $\xi  = x_{0}$ be an arbitrary fixed point in $\Omega $ and $z_0 = Ax_0$.
Then
$$<f-z_{\epsilon }, x_{\epsilon }-x_{0}> = <f-z_{0},
x_{\epsilon }-x_{0}> $$
\begin{equation} \label{p21}
- <z_{\epsilon }-z_{0}, x_{\epsilon }-x_{0}> \leq  <f-z_{0}, x_{\epsilon } -
x_{0}>.
\end{equation}
Taking into consideration a boundedness of $\Omega$ (consequently,
boundedness of the projections $\bar x_{\epsilon }$ and $x_0$),
boundedness of the  maximal
monotone operator $A $ at the point $x_0 \in \;{int} D(A)$
(i.e., the boundedness of $z_0$), we have
from (\ref{p9}), (\ref{p21})
$$||x_{\epsilon } - \bar x_{\epsilon }||^{2} \leq
\epsilon ||f-z_{0}||_{B^*} ||x_{\epsilon }-x_{0}|| = \epsilon
C||x_{\epsilon } - x_{0}||$$
$$\leq  \epsilon C||x_{\epsilon } - \bar x_{\epsilon }|| + \epsilon C
||\bar x_{\epsilon } - x_{0}|| \leq  \epsilon C||x_{\epsilon } - \bar
x_{\epsilon }|| + \epsilon C_{1}, $$
where $C$ and $C_{1}$ are absolute constants. Therefore
$$||x_{\epsilon } - \bar x_{\epsilon }|| \leq  2^{-1}(\epsilon C +
\sqrt { {\epsilon }^2 C^2 + {\epsilon }C_1}),$$
i.e., $||x_{\epsilon }||$  is uniformly bounded.  Then (\ref{p8}) gives
$$<z_{\epsilon } - f, x_{0} - x_{\epsilon }> \ge
\epsilon ^{-1}<J(x_{\epsilon } -  P_{\Omega }x_{\epsilon }), x_{\epsilon }
-x_0> \ge 0 .$$
{}From this
$$<z_{\epsilon }, x_{0} - x_{\epsilon }> \ge
<f, x_{0} - x_{\epsilon }> .$$
Therefore
\begin{equation} \label{p22}
<z_{\epsilon }, x_{\epsilon } - x_{0}> \leq ||f||_{B^*}||x_{\epsilon } -
x_{0}||.
\end{equation}
Lemma  \ref{l2} allows to obtain  the  uniform
boundedness of $z_{\epsilon }$ because
there  exist  the  constant $r_{0} > 0$ and the closed ball
$S(r_{0},x_{0}) \subset  D(F)$ such that for  $x_{\epsilon } \in \Omega$
the estimate
$$<z_{\epsilon }, x_{\epsilon } - x_{0}>  \ge
r_{0}||z_{\epsilon }||_{B^*} - c_{0} (||x_{\epsilon } - x_{0}|| + r_{0}), $$
$$c_{0} = \sup ||A(\xi)||_{B^*} < \infty ,\;\;\; \xi  \in
S(r_{0},x_{0}) ,$$
is satisfied. Consequently, $Ax_{\epsilon }$ are uniformly bounded.
Moreover, it has no more than linear growth on these solutions, i.e.,
$$||Ax_{\epsilon }||_{B^*} \leq r^{-1}_0 (c_{0} + ||f||_{B^*})
||x_{\epsilon } - x_0|| + c_{0}.$$
Then  (\ref{p9}) with $\xi  = \bar x_{\epsilon }$ gives  the  estimate
$||x_{\epsilon } - \bar x_{\epsilon }|| = O(\epsilon )$.

The  sequence $\{x_{\epsilon }\}$  is  weakly  compact,  therefore  one  can
se
subsequence (we do not change  its  notation),  converging weakly  to  some
element $x^*$. The proof of $J(x^* - \bar x^*) = 0 $ (i.e., $x^* \in \Omega $)
is the same as in \cite{l}.

Further, it is easy to see from (\ref{p6}) and  (\ref{p8})  that
$<z_{\epsilon } - f, x - x_{\epsilon }> \ge  0,
\forall x \in \Omega $.  Using the monotonicity of operator $ A$  we
obtain $ (z - f, x - x_{\epsilon }) \ge  0, \forall z \in Ax.$
The weak convergence $x_{\epsilon } \rightarrow  x^*$ gives
$<z - f, x - x^*> \ge  0.$  This inequality  shows
that $x^*$ is a solution of the variational  inequality  (\ref{p1}),
because  it  is
equivalent to (\ref{p2}) \cite{ar}.  Let us now remember the  operator  $A$
is  strictly monotone.  Therefore $x^*$ is a unique solution of (\ref{p1}),
and all  sequences $x_{\epsilon }$
converge weakly to $x^*$.  The theorem is proved.

\begin{re}\label{r5} If $D(A)$ is a bounded set then
the proof of a coerciveness of the operator $T$ can be omitted.
\end{re}
\begin{co}\label{c1}  If the operator $A$  is strictly monotone and
$\Omega$-coercive in $D(A)$, and the set $\Omega $  is  unbounded,
then the assertion  of Theorem  \ref{t1}  is valid.
\end{co}
{\bf Proof.} Let us show that the operator $T$ is $\Omega$-coercive in
$D(A).$  Suppose, that $x_0 \in \Omega$ is the point which describes
the property of $\Omega $-coerciveness of the operator $A.$  Then we have
for all $x \in D(A)$
$$\frac { <Tx, x - x_{0}>}{||x|| } = \frac { <Ax, x - x_{0}>}{||x|| } $$
$$ + \frac { \epsilon ^{-1} <J(x - P_{\Omega }x),
x - x_0>}{||x|| } \ge \frac { <Ax, x - x_{0}>}{||x|| }
\rightarrow \infty ,$$
because
$$<J(x - P_{\Omega }x), x - x_0> \ge  0.$$
The existence theorem for a generalized solution of the equation (\ref{p7})
with the monotone $\Omega $-coercive operator
is proved similarly to \cite {b,ar}.

A uniform boundedness of the solutions $x_{\epsilon }$ follows
from $\Omega $-coerciveness of the operator $A.$ Indeed, if $x_0 \in \Omega$
then (\ref{p8}) gives
$$\frac { <z_{\epsilon }, x_{\epsilon } - x_{0}>}{||x_{\epsilon }|| }
\leq \frac { ||f||_{B^*}||x_{\epsilon } - x_{0}||}{||x_{\epsilon }|| }
\leq ||f||_{B^*} (1 + \frac { ||x_0||}{||x_{\epsilon }|| }) .$$
The assumption $||x_{\epsilon }|| \rightarrow \infty $ contradicts with this
inequality, since its left hand part is bounded as $x_{\epsilon } \rightarrow
\infty $.
Therefore, $x_{\epsilon }$ are bounded  and
$||z_{\epsilon }||_{B^*} \leq C $ from Lemma \ref{l2}.  The further proof
coincides with the proof of Theorem \ref{t1}.

The particular case of this Corollary is the following

\begin{co} \label{c3} If the operator $A$  is strictly monotone and
coercive in $D(A)$, the set $\Omega $  is (possibly) unbounded and
$\theta \in \Omega$, then the assertion  of Theorem  \ref{t1}  is valid.
\end{co}

The assumption of the coerciveness is a strong property of a
structure of the operator $A.$ It permits us to get out of the boundedness
of $\Omega$. Another strong property -  the uniform
monotonicity - leads to strong convergence of the penalty method on the
unbounded sets too.

\begin{co} \label{c2}  If the operator $A$ is uniformly  monotone,
$\psi (t)/t \rightarrow \infty $  as $t \rightarrow \infty $, and $\Omega $
is   unbounded, then  the  solution $x^*$  of  the  variational
inequality (\ref{p1}) exists and the sequence $x_{\epsilon }$ converges
strongly to $x^*$ as $\epsilon \rightarrow  0$.
\end{co}
{\bf Proof.} Using the inequality (\ref{p8}) we obtain for all $x_0 \in
\Omega$
$$<Ax_0 - z_{\epsilon }, x_0 - x_{\epsilon }> \leq
<Ax_0 - f, x_0 - x_{\epsilon }> \leq ||Ax_0 - f||_{B^*}
||x_0 - x_{\epsilon }|| .$$
Then the uniform monotonicity and the local boundedness of $A$ at the
point $x_0$ give
$$ \psi (||x_{\epsilon } - x_0||) \leq C||x_{\epsilon } - x_0|| $$
where $C$ is the constant.

The boundedness of the solution  $x_{\epsilon }$
follows now from the condition
$\psi (t)/t \rightarrow \infty $  as $t \rightarrow \infty .$
Therefore, similarly to Theorem \ref{t1}, the solution $x^*$ exists  and
$ x_{\epsilon } \rightarrow x^* $
weakly. The strong convergence is obtained from the inequality
$$\psi (||x_{\epsilon } - x^*||) \leq <f, x_{\epsilon } - x^*> -
<Ax^*, x_{\epsilon } - x^*> .$$
\\
Theorem \ref{t1} and Corollaries \ref{c1}, \ref{c3} and  \ref{c2} are
generalizations  of  Lions' result \cite{l} (see also \cite{glt}) devoted to
penalty method for variational inequalities with smooth and  bounded
operators.

\section{Variational  inequalities on the approximately given sets.
Stability of the penalty method}
\setcounter{equation}{0}

We  now  formulate  the  results  concerning  the  penalty  method  on
approximately given sets $\Omega _{\sigma } \subset  int\; D(A)$ under
the assumption  that $\Omega _{\sigma }$  are
convex and closed, and Hausdorff distance between $\Omega $  and
$\Omega _{\sigma }$ less than $\sigma,$  i.e.,$\; {\cal H}(\Omega ,
\Omega _{\sigma }) \leq \sigma ,\; 0 \leq  \sigma  \leq  \bar \sigma .$

First of all, we state the lemma about a proximity of  the  generalized
solutions of two penalty equations.  Suppose, that the respective  solutions
$x^{1}_{\epsilon }$ and $x^{2}_{\epsilon }$ of the equations
$$ Ax + \epsilon ^{-1} J(x - P_{\Omega _{1}}x) = f,\;\; \epsilon  > 0 ,$$
$$ Ax + \epsilon ^{-1} J(x - P_{\Omega _{2}}x) = f,\;\; \epsilon  > 0 ,$$
exist and ${\cal H}(\Omega _{1},\Omega _{2}) \leq  \sigma .$

\begin{lem}\label{l4}
Let $B$ be a uniformly  convex  and  uniformly  smooth  Banach
space and let $\delta _{B}(\epsilon)$ be a modulus of convexity of the space
$B$.  Suppose also, that the operator $A$  is uniformly monotone in the form
$$<z_{1}-z_{2}, x_{1}-x_{2}> \ge  \psi (||x_{1}-x_{2}||)||x_{1}-x_{2}||,\;\;
x_1, x_{2} \in  B, \forall z_{1} \in  Ax, \forall z_{2} \in  Ax_{2} $$
where the function $\psi (t)$ is  continuous
on the interval [0, $\infty$), positive for all $t >0$ and
$\psi (0) = 0.$  Then the estimate
\begin{equation} \label{p10}
\psi (||x^{1}_{\epsilon } - x^{2}_{\epsilon }||) \leq  \epsilon ^{-1}C
g^{-1}_{B^*}(2C^{2}L \delta ^{-1}_{B}(4L(d+r)\sigma ))
\end{equation}
occurs,  where $g_{B}(\epsilon ) = \delta _{B}(\epsilon) / \epsilon, \;
g^{-1}_{B}(\cdot)$ is  it's   inverse   function,
$C = 2 {max} \lbrace 1, r+d \rbrace,{\;} 1 < L < 3.18,\; d = {max} \lbrace
d_{1}, d_{2} \rbrace,{\;} d_{i} = {dist} \lbrace \theta ,\Omega _{i} \rbrace,
i = 1,2, $ and \\
$r = \max \lbrace ||x^{1}_{\epsilon }||, ||x^{2}_{\epsilon }||\rbrace. $
\end{lem}
{\bf Proof.} From Definition of the generalized solution it follows
$$<z^1_{\epsilon } - f,x^2_{\epsilon } - x^1_{\epsilon }> +
\epsilon ^{-1} <J(x^1_{\epsilon } - P_{\Omega _{1}}x^1_{\epsilon }),
x^2_{\epsilon } - x^1_{\epsilon }> \ge 0, {\qquad } \forall z^1_{\epsilon }
\in Ax^1_{\epsilon } ,$$
$$<z^2_{\epsilon } - f,x^1_{\epsilon } - x^2_{\epsilon }> +
\epsilon ^{-1} <J(x^2_{\epsilon } - P_{\Omega _{2}}x^2_{\epsilon }),
x^1_{\epsilon } - x^2_{\epsilon }> \ge 0, {\qquad } \forall z^2_{\epsilon }
\in Ax^2_{\epsilon } .$$
Then
$$<z^1_{\epsilon } - z^2_{\epsilon }, x^2_{\epsilon } - x^1_{\epsilon }>
+ \epsilon ^{-1} <J(x^1_{\epsilon } - P_{\Omega _{1}}x^1_{\epsilon }) -
J(x^2_{\epsilon } - P_{\Omega _{2}}x^2_{\epsilon }), x^2_{\epsilon }
- x^1_{\epsilon }> \ge 0 .$$
It is obvious that
$$<J(x^1_{\epsilon } - P_{\Omega _{1}}x^1_{\epsilon }) - J(x^2_{\epsilon }
- P_{\Omega _{2}}x^2_{\epsilon }), x^2_{\epsilon } - x^1_{\epsilon }>$$
$$ = <J(x^1_{\epsilon } - P_{\Omega _{1}}x^1_{\epsilon }) - J(x^1_{\epsilon }
- P_{\Omega _{2}}x^1_{\epsilon }) +
J(x^1_{\epsilon } - P_{\Omega _{2}}x^1_{\epsilon }) - J(x^2_{\epsilon } -
P_{\Omega _{2}}x^2_{\epsilon }), x^2_{\epsilon } - x^1_{\epsilon }>$$
$$ = <J(x^1_{\epsilon } - P_{\Omega _{1}}x^1_{\epsilon }) -
J(x^1_{\epsilon } - P_{\Omega _{2}}x^1_{\epsilon }), x^2_{\epsilon } -
x^1_{\epsilon }>$$
$$ + <J(x^1_{\epsilon } - P_{\Omega _{2}}x^1_{\epsilon }) -
J(x^2_{\epsilon } - P_{\Omega _{2}}x^2_{\epsilon }), x^2_{\epsilon } -
x^1_{\epsilon }>.$$
{}From this one can obtain
$$\epsilon ^{-1}<J(x^1_{\epsilon } - P_{\Omega _{1}}x^1_{\epsilon }) -
J(x^1_{\epsilon } - P_{\Omega _{2}}x^1_{\epsilon }),
x^2_{\epsilon } - x^1_{\epsilon }>  \ge <z^1_{\epsilon } -
z^2_{\epsilon }, x^1_{\epsilon } - x^2_{\epsilon }> $$
$$ \ge  \psi (||x^1_{\epsilon } - x^2_{\epsilon }||)||x^1_{\epsilon } -
x^2_{\epsilon }|| $$
because the operator $\mu (x)$ is monotone, i.e.,
$$<J(x^1_{\epsilon } - P_{\Omega _{2}}x^1_{\epsilon }) -
J(x^2 - P_{\Omega _{2}}x^2_{\epsilon }), x^2_{\epsilon } - x^1_{\epsilon }>
\leq 0 .$$
Now we have
\begin{equation} \label{p22}
 \psi (||x^1_{\epsilon } - x^2_{\epsilon }||) \leq \epsilon ^{-1}
 ||J(x^1_{\epsilon } - P_{\Omega _{1}}x^1_{\epsilon }) -
J(x^1_{\epsilon } - P_{\Omega _{2}}x^1_{\epsilon })||_{B^{*}}.
\end{equation}
This inequality can be continued if to apply (\ref{p18}) to
the right hand part  of (\ref{p22}). Indeed,
$$ ||J(x^1_{\epsilon } - P_{\Omega _{1}}x^1_{\epsilon }) -
J(x^1_{\epsilon } - P_{\Omega _{2}}x^1_{\epsilon })||_{B^{*}} \leq
C_1g^{-1}_{B^*}(2LC_1||P_{\Omega _{1}}x^1_{\epsilon } -
P_{\Omega _{2}}x^1_{\epsilon }||) $$
where $C_1 = \max \lbrace 1, r \rbrace.$   Then
$$\psi (||x^1_{\epsilon } - x^2_{\epsilon }||) \leq \epsilon ^{-1}
C_1g^{-1}_{B^*}(2LC_1||P_{\Omega _{1}}x^1_{\epsilon } -
P_{\Omega _{2}}x^1_{\epsilon }||) .$$
Now Lemma \ref{l3} is to be applied. We have
$$||P_{\Omega _{1}}x^1_{\epsilon } - P_{\Omega _{2}}x^1_{\epsilon }||
\leq  C_2\delta ^{-1}_{B}(4LC_{3}\sigma) ,$$
$$C_2 = 2 \max  \lbrace 1,||x^1_{\epsilon } - P_{\Omega _{1}}x^1_{\epsilon }||,
||x^1_{\epsilon } - P_{\Omega _{2}}x^1_{\epsilon }|| \rbrace, $$
$$C_3 = 2 \max \lbrace ||x^1_{\epsilon } - P_{\Omega _{1}}x^1_{\epsilon }||,
||x^1_{\epsilon } - P_{\Omega _{2}}x^1_{\epsilon }|| \rbrace. $$
Thus,
$$\psi (||x^1_{\epsilon } - x^2_{\epsilon }||) \leq \epsilon ^{-1}
C_1g^{-1}_{B^*}(2LC_1C_2\delta ^{-1}_{B}(4LC_{3}\sigma)).$$
Since
$$||x^1_{\epsilon } - P_{\Omega _{1}}x^1_{\epsilon }||
\leq ||x^1_{\epsilon }||
+ ||P_{\Omega _{1}}\theta || \leq  r + d $$
and
$$||x^1_{\epsilon } - P_{\Omega _{2}}x^1_{\epsilon }||
\leq ||x^1_{\epsilon }|| +
||P_{\Omega _{2}}\theta || \leq  r + d, $$
the estimate (\ref{p10}) holds. Lemma is proved.

\begin{re} \label{r5}
If $\psi (t)$ is a continuous, strictly increasing function,
$\psi (0) = 0,$  then the estimate
$$ ||x^{1}_{\epsilon } - x^{2}_{\epsilon }|| \leq \psi ^{-1}[ \epsilon ^{-1}C
g^{-1}_{B^*}(2C^{2}L \delta ^{-1}_{B}(4L(d+r)\sigma ))] $$
is valid.
\end{re}

\begin{re} \label{r2}
In \cite{a1} we noted that the  estimate  of  Hausdorff  distance
${\cal H}(\Omega ,\Omega _{\sigma }) < \sigma $ has its own importance,
mainly, for bounded  sets.   Therefore,
new  conditions,  characterizing  the  measure   of
proximity of two (possibly unbounded) sets, have  been  introduced:
\begin{equation} \label{p11}
d(x,\Omega ) \leq  \sigma f_{1}(||x||),\;\; \forall x \in  \Omega _{\sigma },
\end{equation}
\begin{equation} \label{p12}
d(x,\Omega _{\sigma }) \leq  \sigma f_{2}(||x||),\;\; \forall x \in  \Omega,
\end{equation}
where $f_{1}(t)$ and $f_{2}(t)$ are finite functions for all $t \ge  0,
d(x,\Omega ) = \inf_{y \in \Omega} ||x-y||.$
If (\ref{p11}) and (\ref{p12}) are true, then in formula  (\ref{p10})
$d+r$  under  the  sign  of
function $\delta ^{-1}_{B}(\cdot )$ is replaced by the value $(d+r)
(f_{1}(2r+d)+f_{2}(2r+d))$.
\end{re}

\begin{re} \label{r11}
A correspondence between the modulus of a convexity of the space $B^*$
and  the modulus of a smoothness of the space $B$ can be given in a form
\cite{d}:
$$\rho _B (\tau) \ge \epsilon \rho /2 - \delta _{B^{*}} (\epsilon),\;\;
0 \leq \epsilon \leq 2,\;\; \tau > 0 .$$
{}From this we obtain the following estimate
$$||Jx-Jy||_{B^{*}} \leq 4C h_B (8CL||x - y||),$$
where $h _B (\tau) = \rho _B (\tau) /\tau $ and
$$ C = C(||x||,||y||) = 2 {max} \lbrace 1, \sqrt{(||x||^2 + ||y||^2)/2}
\rbrace .$$
Then, under the conditions and notations of Lemma \ref{l4} we have
$$\psi (||x^{1}_{\epsilon } - x^{2}_{\epsilon }||) \leq  \epsilon ^{-1}C
h^{-1}_{B}(2C^{2}L \delta ^{-1}_{B}(4L(d+r)\sigma )) .$$
\end{re}

Let us consider the penalty equation on the sets $\Omega _{\sigma }$
\begin{equation} \label{p13}
Ax + \epsilon ^{-1}J(x-P_{\Omega _{\sigma }}) = f,\;\; x \in  D(A),\;
\epsilon > 0,\; \sigma \ge 0.
\end{equation}
\begin{th} \label{t2}
Let the assumptions of Lemma \ref{l4} be satisfied, provided that either
${\cal H}(\Omega , \Omega _{\sigma })
\leq \sigma $ holds or (\ref{p11}), (\ref{p12}) occur,
and $\psi (t) \rightarrow  \infty $ as $t \rightarrow \infty $.
Then:\\
(i) the solution $x^*$ of variational inequality (\ref{p1}) exists; \\
(ii)  the generalized solutions $x^{\sigma }_{\epsilon }$ of the equation
(\ref {p13}) exist for all $\sigma  \ge  0,\; \epsilon > 0$; \\
(iii) the condition
\begin{equation} \label{p14}
g^{-1}_{B^*}(\delta ^{-1}_{B}(\sigma))/\epsilon \rightarrow 0
\end{equation}
involves a strong convergence of the penalty method, i.e.,
$$||x^{\sigma }_{\epsilon } - x^*||_{B} \rightarrow 0\;\;{as\;}
\epsilon \rightarrow 0,\; \sigma \rightarrow 0.$$
\end{th}

{\bf Proof.} An existence of the soluion $x^*$ of variational inequality
(\ref{p1}) and the generalized solutions $x^{\sigma }_{\epsilon }$ of the
equation (\ref {p13}) are guaranteed by Corollary \ref{c1}, because
$\psi (t) \rightarrow  \infty $ as $t \rightarrow \infty $.

Let  $x_{\epsilon }$  is a generalized soluion  of the  equation (\ref{p7}).
Then
$$||x^{\sigma }_{\epsilon } - x^*|| \leq ||x^{\sigma }_{\epsilon } -
x_{\epsilon } + x_{\epsilon } - x^*||
\leq ||x^{\sigma }_{\epsilon } - x_{\epsilon }|| + ||x_{\epsilon } - x^*||.$$
Note that the
properties of the functions $\psi (t), \delta (\epsilon ),
\delta ^{-1}_{B}(\cdot )$ and $g^{-1}_{B^*}(\cdot )$ such that
$\psi (t) \rightarrow  0 $ as $t \rightarrow 0 $ and
$g^{-1}_{B^*}(\delta ^{-1}_{B}(\sigma)) \rightarrow 0 $ as $ \sigma
\rightarrow 0.$ Therefore
$||x^{\sigma }_{\epsilon } - x_{\epsilon }|| \rightarrow 0$ as $ \epsilon
\rightarrow 0,\; \sigma \rightarrow 0$ and (\ref{p14}) is satisfied
(see Lemma \ref{l4}). Furthermore
$||x_{\epsilon } - x^*|| \rightarrow 0$ as $ \epsilon \rightarrow 0$
(see Corollary \ref{c1}). The theorem is proved.\\
\\

Let us suppose further that the operator  $A$  and  the  element $f$  are
perturbed, such that the penalty equation
\begin{equation} \label{p15}
A^{h}x + \epsilon ^{-1}J(x-P_{\Omega _{\sigma }}x ) = f^{\omega },\;
\epsilon  > 0,\; \epsilon \rightarrow 0
\end{equation}
is given under the additional conditions:\\
\\
(a)$\;$ $A^{h}$ are maximal monotone operators too,$\; 0 \leq  h \leq
\bar h,\;
D(A^{h}) = D(A),\; \theta  \in  D(A), \\
{\cal H}_{B^*}(M^{h}x,Mx) \leq  h\gamma (||x||)$,  where ${\cal H}_{B^*}
(Q_{1},Q_{2})$  is  the Hausdorff  distance between the sets
$Q_{1}$ and $Q_{2}$ in the space $B^*, Mx$ and $M^{h}x$ are the  ranges
of the operators $A$  and $A^{h}$ at the points $x \in
\bigcup {}_{\sigma \ge 0}\Omega _{\sigma }$,  and $\gamma (t)$  is  a
continuous nondecreasing function for $t \ge  0$;\\
\\
(b)$\;$ $||f^{\omega } - f||_{B^{*}} \leq \omega ,\; 0 \leq  \omega
\leq \bar \omega . $

\begin{th} \label{t3}
Under  conditions  (a)  and  (b)  and  the  hypothesis  of
Theorem \ref{t2}, the generalized solutions $x^{\Delta }_{\epsilon }$
of  equation  (\ref {p15})  exist  and
uniformly bounded for all $\Delta  = (\sigma ,h,\omega ),\; \epsilon  > 0,$
and $||x^{\Delta }_{\epsilon } - x^*||_{B} \rightarrow 0 $ as $\Delta
\rightarrow 0,\; \epsilon \rightarrow 0. $
\end{th}
{\bf Proof.}
For generalized solutions $x_{\epsilon }$ and $x^{\Delta }_{\epsilon }$
of the equations  (\ref {p7})  and (\ref {p15}) we have
$$\psi (||x^{\Delta }_{\epsilon } - x_{\epsilon}||) \leq  \epsilon ^{-1}
||J(x^{\Delta }_{\epsilon } - P_{\Omega _{\sigma }} x^{\Delta} _{\epsilon} -
J(x^{\Delta }_{\epsilon } - P_{\Omega } x^{\Delta} _{\epsilon })||_{B^{*}} +
\omega  +  h\gamma (||x_{\epsilon}||).$$
Here we used the following relations: $\forall z_{\epsilon } \in A
x_{\epsilon }, \forall z^{\Delta }_{\epsilon } \in A
x^{\Delta }_{\epsilon },
\forall z^{\Delta ,h}_{\epsilon } \in A^{h}x^{\Delta }_{\epsilon },$

$$<z^{\Delta ,h}_{\epsilon } - z_{\epsilon }, x^{\Delta }_{\epsilon } -
x_{\epsilon }> =  <z^{\Delta }_{\epsilon } - z_{\epsilon },
x^{\Delta }_{\epsilon } - x_{\epsilon }> +
<z^{\Delta ,h}_{\epsilon } - z^{\Delta }_{\epsilon },
x^{\Delta }_{\epsilon } - x_{\epsilon }>$$
$$\ge  \psi (||x^{\Delta }_{\epsilon } - x_{\epsilon }||) ||x^{\Delta }
_{\epsilon } - x_{\epsilon }|| + <z^{\Delta ,h}_{\epsilon } - z^{\Delta }
_{\epsilon }, x^{\Delta }_{\epsilon } - x_{\epsilon }> , $$

$$<z^{\Delta ,h}_{\epsilon } - z^{\Delta }_{\epsilon }, x^{\Delta }
_{\epsilon } - x_{\epsilon }> \leq  h\gamma (||x^{\Delta }_{\epsilon }||)
||x^{\Delta }_{\epsilon } - x_{\epsilon }||,$$
$$<f^{\omega } - f,
x^{\Delta }_{\epsilon } - x_{\epsilon }> \leq \omega
||x^{\Delta }_{\epsilon } - x_{\epsilon }||, $$
and
$$<J(x_{\epsilon } - P_{\Omega }x_{\epsilon }) - J(x^{\Delta }_{\epsilon } -
P_{\Omega }x^{\Delta }_{\epsilon }), x_{\epsilon } - x^{\Delta }_{\epsilon }>
\ge  0 .$$
Therefore, $||x^{\Delta }_{\epsilon } - x_{\epsilon }|| \rightarrow 0,$
because $x_{\epsilon }$ and $x^{\Delta }_{\epsilon }$ are uniformly bounded,
$h \rightarrow 0, \omega \rightarrow 0,$ and (\ref{p14}) is fulfilled.
We finish the proof by adding Corollary \ref{c2}. \\
\\
Finally,  consider  the  variational  inequality  (\ref{p1}),  provided   its
solution set is nonempty and $A$ is  a  strictly  proper  monotone  operator.
Then (\ref{p1}) belongs to the class of unstable nonlinear  problems.
Therefore, in the penalty method it is necessary to  apply  some
regularization.   In this connection we use a duality mapping
$Jx$ with vanishing parameter $\alpha $  as the regularizing operator,
\begin{equation} \label{p16}
A^{h}x + \epsilon ^{-1} J(x - P_{\Omega _{\sigma }}x) + \alpha Jx =
f^{\omega },\;\; \epsilon  > 0,\;\; \alpha  > 0.
\end{equation}
\begin{th} \label{t3}
Let $B$  be a uniformly convex and  uniformly  smooth  Banach
space, $\delta _{B}(\epsilon )$  be a modulus of convexity of B. Let either
${\cal H}(\Omega ,\Omega _{\sigma }) \leq  \sigma $   or
(\ref{p11}), (\ref{p12}) occur.  Then:\\
(i)$\;$ the generalized solutions $x^{\Delta }_{\alpha ,\epsilon }$ of equation
(\ref{p16}) exist for all $\epsilon  > 0,\;
\alpha  > 0,$  and they are uniformly bounded;\\
(ii) $\;x^{\Delta }_{\alpha ,\epsilon }$  converge strongly to $x^*$  under
the conditions
\begin{equation} \label{p17}
{\frac {h + \omega}{\alpha}} \rightarrow 0,\;\;
{\frac {g^{-1}_{B^*}(\delta ^{-1}_{B}
(\sigma))}{\alpha \epsilon}} \rightarrow 0,\;\; \alpha \rightarrow 0,\;\;
\epsilon  \rightarrow 0.
\end{equation}
\end{th}
{\bf Proof.} For the difference of two  solutions $x^{0}_{\alpha ,\epsilon }$
and $x^{\Delta }_{\alpha ,\epsilon }$  of  the
equation (\ref{p16}) with $A^{h} = A, f^{\omega } = f, \Omega _{\sigma } =
\Omega $ and $A^{h},f^{\omega }$ and $\Omega _{\sigma }$, respectively, we
have
$$||x^{0}_{\alpha ,\epsilon } - x^{\Delta }_{\alpha ,\epsilon }|| \leq
C_{2}g^{-1}_{B}[2C_{2}L(\alpha ^{-1}(h\gamma (r) + \omega ) +
(\alpha \epsilon )^{-1}Cg^{-1}_{B^*}(2C^{2}L \delta ^{-1}_{B}(4L(d+r)\sigma
)))], $$
where $C = 2 \max \lbrace 1,r+d \rbrace ,\; C_{2} = 2\max \lbrace 1,r \rbrace,
r = \max \lbrace ||x^{0}_{\alpha ,\epsilon }||,
||x^{\Delta }_{\alpha ,\epsilon }|| \rbrace ,
{\cal H}(\Omega ,\Omega _{\sigma }) \leq \sigma .$

Therefore, $||x^{\Delta }_{\alpha ,\epsilon } - x^{0}_{\alpha ,\epsilon }||
\rightarrow 0 $  by virtue of the first two relations in  (\ref{p17}),
while $\alpha \rightarrow 0 $ and $\epsilon \rightarrow 0 $ guarantee  the
convergence $x^{0}_{\alpha ,\epsilon }$  to $x^* $ \cite{a2}.   The
boundedness of $||x^{0}_{\alpha ,\epsilon }||$ and
$||x^{\Delta }_{\alpha ,\epsilon }||$
is proved by analogy with Theorem  \ref{t1}.

We considered here the condition ${\cal H}(\Omega ,\Omega _{\sigma }) \leq
\sigma $; for (\ref{p11}) and (\ref{p12})  one  can  see
Remark \ref{r2}.

\begin{re} \label{r6}
If the operator $A$ is an arbitrary monotone, then
in corresponding statements above, we suppose that its
maximal monotone expansion in $D(A)$ is coercive and uniformly  monotone.
\end{re}

\begin{re} \label{r7}
The assumption $\Omega  \subseteq  {int} D(A), \; {int} \Omega  \neq
\emptyset $ can be replaced with  $\Omega  \subseteq  D(A)$ and
$ {int}D(A) \bigcap \Omega \neq \emptyset.$
\end{re}

\begin{re} \label{r4}
All the results presented above relate to the  problems  of
conditional and unconditional optimization for convex nonsmooth functionals
with (possibly) unbounded gradients or subgradients.
Lions  notes in \cite{l} that there are monotone continuous coercive and
unbounded operators. For example, a gradient of the functional
$\psi = \sum _{n=1}^{\infty} |u_m|^{m+1}/(m+1)$ in the space $l^2$ of the
sequences $ \lbrace u_m| m= 1,2,... \rbrace $  is not bounded on the  bounded
sequences  $f_1 = \lbrace 2,0,0,... \rbrace, f_2 = \lbrace 0,2,0,0,...
\rbrace $ etc.
\end{re}

Multiple applications of the penalty methods for the  solving variational
inequalities with differential operators can be found in \cite{l,ks,glt}.


\begin{thebibliography}{99}

\bibitem{a1}  Ya.I. Alber, The regularization method  for  variational
inequalities with nonsmooth unbounded operators in Banach space,
{\it Appl. Math. Lett., \/} 6 (1993), 63-68.

\bibitem{a2}  Ya.I. Alber,  The  solution  of  nonlinear  equations  with
monotone
operators in Banach spaces, {\it Siberian Math. J., \/}  16 (1975), 1-8.

\bibitem{an3}  Ya.I. Alber, A.I. Notik, Perturbed unstable variational
inequalities with unbounded  operators  on  approximately  given  sets,
{\it Set-Valued Analysis \/} (accepted for publication).

\bibitem{an1} Ya.I. Alber and A.I. Notik, Geometric properties of
Banach spaces and approximate
methods for solving nonlinear operator equations, {\it Soviet Math.
Dokl.\/,} 29 (1984), 611-615 .

\bibitem{an2}  Ya.I. Alber, A.I. Notik, Parallelogram inequalities in Banach
spaces and some properties of a duality  mapping, {\it Ukrainian Math. J. \/,}
40 (1988), 650-652 .

\bibitem{an4} Ya.I. Alber and A.I. Notik, On some estimates for
projection operator in Banach space (to appear).

\bibitem{ar}  Ya.I.  Alber,  I.P.  Pjazantseva,   Variational   inequalities
with discontinuous  monotone  mapping,  {\it Soviet  Math.  Dokl., \/}
25 (1982), 206-210.

\bibitem{b}   F.E. Browder,  Nonlinear maximal monotone operators in Banach
space,  {\it Math. Annalen, \/} 175 (1968), 89-113.

\bibitem{d}  D. Distel', The Geometry of Banach Spaces, {\it Lecture Notes
Math., \/}  485, Springer, 1975.

\bibitem{fi}  T. Figiel,  On the moduli of convexity and smoothness,
{\it Studia Mathematica, \/} 56 (1976), 121-155.

\bibitem{glt} R.G. Glowinski, J.-L. Lions and R. Tremolieres, {\it Alalyse
Numerique des Inequations Variationnelles, \/} Vol.1, Dunod, Paris, 1976.

\bibitem{ks} D. Kinderlehrer and G. Stampacchia, {\it An  Introduction to
Variat
Inequalities and Their Applications, \/} Academic Press, New York, London,
Toronto, 1980.

\bibitem{l} J.-L. Lions,
{\it Quelques methodes de resolution des problemes aux limites non
lineaires,\/}  Dunod Gauthier-Villars, Paris, 1969.

\bibitem{p} G. Pisier, Martingales with values in uniformly convex spaces,
{\it Isr. J. Math., \/} 20 (1975), 326-350.

\bibitem{rc}  R.T. Rockafellar, Local boundedness of nonlinear monotone
operators,  {\it Tran. Amer. Math. Soc., \/} 149 (1970), No. 1, 75-88.

\end{thebibliography}
\end{document}